%% file: RD_FR_Connection_Paper.tex
\documentclass[USenglish]{article}
\usepackage{fullpage}
\usepackage{amsmath,amsfonts,amsthm,amssymb}
\usepackage{color,xcolor}
\usepackage{ifpdf}
\usepackage{psfrag}
\usepackage{graphicx,graphics,subfigure}
\usepackage{url}
\usepackage{hyperref}
\usepackage{todonotes}
\usepackage{ulem} 
\newtheorem{theorem}{Theorem}[section]

\newtheorem{proposition}[theorem]{Proposition}
\newtheorem{remark}[theorem]{Remark}

\newtheorem{ex}[theorem]{Example}

\newcommand{\R}{\mathbb R}

\newcommand{\PP}{\mathbb P}
\newcommand{\TT}{\mathcal{T}}

\newcommand{\EE}{\mathcal{E}}

\newcommand{\ee}{e}

\newcommand{\vecn}{\boldsymbol{\mathbf n}}

\newcommand{\E}{\mathrm{E}}

\newcommand{\bbf}{{\mathbf {f}}}
\newcommand{\bxx}{{\mathbf {x}}}
\newcommand{\bn}{{\mathbf {n}}}

\newcommand{\Bf}{\mathbf{f}}

\newcommand{\bu}{\mathbf{u}}
\newcommand{\bv}{\mathbf{v}}

\newcommand{\bbv}{\mathbf{v}}

\newcommand{\bbu}{\mathbf{u}}
\newcommand{\bbg}{\mathbf{g}}
\newcommand{\hbbg}{\hat{\mathbf{g}}}
\newcommand{\hbbf}{\hat{\mathbf{f}}}
\newcommand{\hbbfe}{\hat{\mathbf{f}}_{\sigma,\sigma'}}

\newcommand{\red}[1]{\textcolor{red}{#1}}

\newcommand{\bs}{\boldsymbol}
\renewcommand{\div}{\operatorname{div}}
\newcommand{\est}[1]{\left\langle#1\right\rangle}
\newcommand{\mean}[1]{\overline{#1}}
\newcommand{\bbfh}{{\mathbf {f}^h}}
\newcommand{\Ol}{\mathcal{O}}

\hypersetup{
pdftitle={On the Connection between Residual Distribuation Schemes and the Flux Reconstruction Approach.\\
Application to entropy conservative and entropy dissipative schemes}
pdfauthor={P. \"Offner},
pdfpagemode=UseOutlines,
linkbordercolor=0 0 0,
linkcolor=red,
citecolor=blue,
colorlinks=true,
bookmarks = true
}
\begin{document}
\title{On the Connection between Residual Distribution Schemes and Flux Reconstruction 
}
\author{R. Abgrall, E. le M\'el\'edo, and P. \"Offner\\
Institute of Mathematics,
University of Zurich, Switzerland}
\date{
\today}
\maketitle

\begin{abstract}
In this short paper, we are considering the connection between
the \emph{Residual Distribution Schemes} (RD) and the \emph{Flux Reconstruction}
(FR) approach.
We demonstrate that  flux reconstruction can be recast into the RD
framework and vice versa. Because of this close connection we are able to apply known results from RD schemes to FR methods.
In this context we 
propose a first demonstration of
entropy stability for the FR schemes under consideration and 
show how to construct entropy stable numerical schemes based on our FR methods.
Simultaneously, we do not restrict the mesh to  tensor 
structures or triangle elements, but rather allow polygons. 
The key of our analysis is a proper 
choice of the correction functions for which we present an approach here.
\end{abstract}
\tableofcontents
\input{1_Introduction_Paper}

\input{2_FR_to_RD}

\input{3_Transformation}

\input{4_Summary}

\input{Appendix}

\section*{Acknowledgements}
The second and third authors have been funded in by the the SNF project (Number 175784) ``
Solving advection dominated problems with high order 
schemes with polygonal meshes: application to compressible and incompressible flow problems''.

\bibliographystyle{abbrv}
\bibliography{literature}

\end{document}

%% file: 1_Introduction_Paper.tex
\section{Introduction}
We are interested in the approximation of mainly non-linear hyperbolic problems 
like Euler equations or the MHD equations.
The construction of high order methods
for these problems is widely studied in current research 
with the aim to find good ways to build preferable schemes. 
All of these methods have in common that they are based of 
either a finite difference (FD) or a finite element (FE) approach.
In the last years, great efforts have been made to transform 
numerical schemes from one to another and to use techniques 
which are originally used in a  different framework.
Here, the summation-by-parts (SBP) operators are a good example to mention.
SBP operators originate in the FD framework \cite{kreiss1974finite}
and lead to an ansatz to prove stability in a  way similar
 to the continuous analysis (see \cite{svard2014review,  fernandez2014review, hicken2016multidimensional} 
 and the references therein). 
In \cite{gassner2013skew} the author transforms the
technique to a Discontinuous Galerkin (DG) 
spectral element method (DGSEM) using the nodes of
Lobatto-Legendre quadrature, 
and in 
\cite{oeffner2018error, oeffner2017stability, ranocha2018generalised, ranocha2018stability, ranocha2016summation, ranocha2017extended} 
SBP operators are applied to the Correction Procedure via Reconstruction (CPR) or 
Flux Reconstruction (FR) method to extend stability proofs in a more general framework.

In this paper, we also deal with the reinterpretation
or transformation of two classes of 
numerical methods,  meaningly the Residual Distribution (RD) schemes and the FR methods. 
Both  lead to a general framework which contains several numerical schemes.
The FR creates a unifying framework for several high-order methods such as 
DG, spectral difference (SD) and spectral volume (SV) methods. 
These connection are already pointed out in \cite{huynh2014high}
and in the references therein.
However, since the early work of Roe, Deconinick and Struijs
\cite{roe1986characteristic, struijs1991fluctuation}, 
the RD schemes have been further developed in a series of papers, 
e.g. \cite{abgrall2001toward, abgrall2003high, deconinck2000status, ricchiuto2007application}.  
A connection between RD to DG is explained in \cite{abgrall2009development} and, 
because of the close relation between RD and FR to DG, 
it seems natural to study the link between RD and FR.

Besides accuracy and robustness of the numerical scheme, another 
desirable property of numerical schemes is entropy stability.
Recently, efforts have been made to construct numerical methods 
enjoying entropy stability. 
So far, linear stability remains mainly investigated in the context of FR, as e.g. 
\cite{castonguay2012newclass, castonguay2013energy, jameson2012nonlinear, vincent2015extended,
wang2016review, vermeire2016properties, williams2014energy}.
By embedding FR into the RD framework we are able to follow
the steps of \cite{abgrall2017general} 
and construct FR schemes that are 
also entropy conservative. The key is a proper choice of the 
correction functions, which is discussed in this paper.
Another beautiful consequence of our abstract approach 
it that we do not need to restrict our mesh in two dimensions on triangles 
or tensor structures. Our approach is valid for general polygons, 
extending  the current results. \\
The paper is organized as follows.
In the next section we shortly repeat the main idea of RD schemes and FR schemes. We 
introduce the notations that will be used later in this work. 
After that we explain  the flux reconstruction approach and formulate the 
schemes in the RD context.
Therefore, the definition of the correction functions in FR are essential and 
two conditions guaranteeing conservation are derived.
Furthermore, in the section \ref{sec:Transformation}, we transform 
theoretical convergence and stability  results from RD to our FR schemes.
Simultaneously, we make some preparations for discussing  entropy stability.
In the next section \ref{sec:entropy}  we follow the steps from \cite{abgrall2017general} 
and construct an entropy conservative/stable numerical scheme based on our 
FR approach. Simultaneously, by bringing our investigations together we are able to derive straight conditions on our correction functions so that the resulting
FR methods are naturally entropy conservative. 
We then summarize everything and conclude on the admissibility
criteria for the correction functions.
In the appendix we further extend the investigation on entropy stability and correction functions. 
Entropy stability is usually associated with the condition of Tadmor on the 
numerical flux. 
Following the study of \cite{abgrall2017general, abgrall2017some, abgrall2017some2} where similar conditions for
RD schemes are proposed and making use of the link between RD and FR,
we are able to derive conditions on FR schemes that guarantee entropy stability.
It is shown that the correction functions have to fulfill an inequality which is derived from an entropy inequality.  
 However, if this more a theoretical condition can be checked, it does not 
lead to a construction method for entropy conservative/stable schemes yet. 
Nevertheless, as it may be interesting for future research, we detail it in the appendix for the sake of completeness.

\section{Residual Distribution Schemes and Flux Reconstruction}

\subsection{Residual Distribution Schemes - Basic Formulation}\label{sec:RD_Basic}

We follow the ideas and notations for the RD methods of 
\cite{abgrall2006residual, abgrall2017general, abgrall2017some, abgrall2017some2, abgrall2017high}.
As already described in \cite{abgrall2017some}, for example, the 
\emph{Discontinuous Galerkin Method} (DG) 
can be interpreted as a RD scheme. Thus, by the
close connection between  DG and FR (see \cite{huynh2007flux}), such an interpretation
is also possible for FR schemes.

In this paper we are considering the steady state problem 
\begin{equation}\label{eq:steady_problem}
 \div \bbf (\bu)= \sum\limits_{j=1}^d \frac{\partial \bbf_j}{\partial x_j}\left(\bu \right)=0 \text { for } \bxx\in \Omega \subset \R^d 
\end{equation}
together with the initial condition 
\begin{equation}\label{eq:initial condition}
 \left( \nabla_{\bu} \bbf (\bu) \cdot \vecn(\bxx) \right)^-\left( \bu- \bu_b \right)=0 \text{ on } \partial \Omega,
\end{equation}
where $\vecn(\bxx)$ is the outward normal vector at
$\bxx\in \partial \Omega$ and $\bu_b $
is a regular enough function.
An extension of \eqref{eq:steady_problem} to unsteady problems is
straightforward by following the steps from 
\cite{abgrall2017high}.
The flux function is given by 
\begin{equation*}
 \bbf_j = \begin{pmatrix}
            f_{1,j}\\
            \vdots\\
            f_{p,j}
           \end{pmatrix} 
           \text{ and } \bu =\begin{pmatrix}
            u_{1}\\
            \vdots\\
            u_{p} 
              \end{pmatrix} \in D\subset \R^p 
\end{equation*}
is the conserved variable.

Later on, we will focus on the entropy $\bs U$, which fulfils  the condition 
\begin{equation}\label{eq:Entrop}
 \nabla_{\bu} \bs U \cdot \nabla_u \bbf_j = \nabla_{\bu}  g_j, \quad \forall j=1,\dots, d
\end{equation}
where $\bs g=(g_1, \cdots,g_d)$, $g_j\in C^1(\Omega)$ and 
is called the entropy flux function.
If $\bu$ is smooth, then the additional conservation relation
\begin{equation}\label{eq:Entrop_div}
 \div \bs g( \bu ) = \sum\limits_{j=1}^d \frac{\partial g_j}{\partial x_j} (\bu)=0
\end{equation}
holds. If $\bu$ is a weak entropy 
solution of \eqref{eq:steady_problem}-\eqref{eq:initial condition}, 
then $\div \bs{g}(\bu) \leq 0$ holds in the sense of distribution.

We split $\Omega$ in a partition of elements $K$, and 
approximate our solution in each element by a polynomial of degree $k$.
The numerical solution is denoted $\bu^h$.
Therefore, the numerical solution lies 
in the space $V^h=\bigoplus\limits_{K}\{\bu^h \in L^2(K)^d, \,\bu^h|_K\in \PP^k(K) \}$. 

In addition, we set $\{ \phi_\sigma \}_{\sigma \in \sum_K}$ a set  of 
basis functions  for the space $\PP^k(K)$, where 
$\sum_K$ is a set of degrees of freedom of linear forms acting on the set $\PP^k(K)$,
which will be used in every element to express $\bu^h$. 
Considering all the elements covering $\Omega$,
the set of degrees of freedom is denoted by $S$ and a generic degree of
freedom by $\sigma$. Furthermore, for any $K$ we have that
\begin{equation*}
 \forall \bxx \in K, \qquad \sum\limits_{\sigma\in K} \phi_\sigma(x)=1.
\end{equation*}
The key of the RD schemes is  to define residuals $\Phi_\sigma^{K}(\bu^h)$
on every element $K$, satisfying element-wise the following conservation relation
\begin{equation}\label{eq:residual_element}
 \sum\limits_{\sigma\in K} \Phi^K_\sigma (\bu^h) = \oint_{\partial K}  \hbbf( \bu^h, \bu^{h,-} )\operatorname{d} \gamma,
\end{equation}
where $\bu^{h,-}$ is the approximated solution on the 
other side of the local edge/face of $K$, $\hbbf$
is a consistent numerical flux, i.e. 
$\hbbf(\bu, \bu)= \bbf( \bu) \cdot \vecn$, and 
$\oint_K$ is the boundary integral evaluated by a numerical quadrature rule. 
Simultaneously, we have to consider the residuals on the boundary elements
$\Gamma$. For any degree of freedom belonging to the boundary $\Gamma$, we assume 
that $\Phi_\sigma^{\Gamma}(\bu^h)$ fulfils the conservation relation 
\begin{equation}\label{eq:residual_boundary}
 \sum\limits_{\sigma\in \Gamma} \Phi^\Gamma_\sigma (\bu^h) = \oint_{\partial \Gamma}\left(   \hbbf( \bu^h, \bu_{b} ) -\bbf( \bu^h)\cdot \vecn \right) \operatorname{d} \gamma.
\end{equation}
The discretisation of \eqref{eq:steady_problem}-\eqref{eq:initial condition} is given by the following formula.
For any $\sigma \in S$, it reads
\begin{equation}\label{eq:Discretisatin_RD}
 \sum\limits_{K\subset \Omega, \; \sigma \in K} \Phi^K_\sigma (\bu^h) +\sum\limits_{ \Gamma \subset \partial \Omega, \; \sigma \in \Gamma} \Phi_\sigma^\Gamma (\bu^h) =0. 
\end{equation}
We are able to embed the discretisation  \eqref{eq:Discretisatin_RD} into several numerical methods like finite element or DG depending on the solution space $V^h$
and the definition of the residuals, see
\cite{abgrall2017some} for details. 
Here, we repeat it shortly for the DG scheme before  we extend it to the FR methods. 

A weak formulation of DG reads:
Find $\bu^h \in V^h$ such that for any $\bs v^h \in V^h$,
\begin{equation}\label{eq:DG_Variational}
 \begin{aligned}
   a(\bu^h, \bs v^h):=&
\sum\limits_{K\subset\Omega}  \left(-\oint_K \nabla \bs v^h \cdot \bbf(\bu^h) \operatorname{d} x +  \oint _{\partial K} \bs v^h \cdot \bs \hbbf (\bu^h, \bu^{h,-} ) \operatorname{d} \gamma \right)
\\
&+ \sum\limits_{\Gamma \subset \partial \Omega} \oint_\Gamma \bs v^h \cdot \left( \bs \hbbf(\bu^h, \bu_b)-\bbf(\bu^h)\cdot \vecn \right) \operatorname{d} \gamma=0.
 \end{aligned}
\end{equation}
 Here we have defined for the boundary faces $\bu^{h,-}=\bu^h$ and used 
 the fact that the expression $ \nabla \bs v^h \cdot \bbf(\bu^h)$ implies
\begin{equation}\label{eq:nabla_not}
 \nabla \bs v^h \cdot \bbf(\bu^h)=\sum\limits_{j=1}^d \left(\frac{\partial }{\partial x_j} \bs v^h  \right) \cdot \bbf_j(\bu^h).
\end{equation}
The strong version of DG is obtained by applying integration-by-parts (summation-by-parts in the discrete sense) another time.
We obtain the corresponding RD scheme's residuals by comparing
\eqref{eq:Discretisatin_RD} and \eqref{eq:DG_Variational}.
For the inner elements, we get
\begin{equation}\label{eq:DG_Residual_inner}
\Phi_\sigma^{K,DG} (\bu^h) = -\oint_K \nabla  \phi_\sigma \cdot \bbf(\bu^h) \operatorname{d} \bxx +  
\oint _{\partial K} \phi_\sigma \cdot \bs \hbbf (\bu^h, \bu^{h,-} ) \operatorname{d} \gamma. 
\end{equation}
The boundary residuals are given by 
\begin{equation}\label{eq:DG_Residual_boundary}
 \Phi_\sigma^{\Gamma,DG}(\bu^h)=\oint_\Gamma   \phi_\sigma \cdot \left( \bs \hbbf(\bu^h, \bu_b)
 -\bbf(\bu^h)\cdot \vecn \right) \operatorname{d} \gamma.
\end{equation}

Note that the expressions \eqref{eq:DG_Residual_inner} and \eqref{eq:DG_Residual_boundary}
satisfy  the conservation relations
\eqref{eq:residual_element} and \eqref{eq:residual_boundary}. 

In view of later use, we  set an  expression for the average of 
the left and right states of $a$ on the boundary of $K$ and a jump condition that
for any function $\omega$ as follows.
\begin{equation}\label{Eq:Jump_Average}
 \{a\}:=\frac{1}{2}\left(a^K+a^{K,-} \right) \qquad [\omega]:=\omega_{|K^-}-\omega_{|K}
\end{equation}

\subsection{Flux Reconstruction Approach on Triangles}\label{subsec:FR}
In our research we focus on two dimensional problems and target the use general polygonal meshes.
Therefore, we start by introducing the FR approach directly  on triangles following the explanations of \cite{castonguay2012newclass}.
For a detailed introduction into FR we strongly recommend the review article \cite{huynh2014high} and the references therein.
Instead of considering our steady state problem \eqref{eq:steady_problem}, we focus in this subsection on a 
two dimensional scalar conservation law 
\begin{equation}\label{eq:hyp_prob}
 \partial_t u+\div( \bbf)=0
\end{equation}
within an arbitrary domain $\Omega$.  The flux is now $\bbf=(f(u),g(u))$ and $\div$ is divergence in the space variables $x,y$. 
The domain is splitted into $N$ non-overlapping elements $\Omega_k$ such that $\Omega=\bigcup_{k=1}^K \Omega_k$. 
We use here a conforming triangulation of $\Omega$. Note that a quadrangulation would also be possible.
Both $u$ and $\bbf$ are approximated by polynomials in every element $\Omega_k$ and  their  total approximation in $\Omega$ is given by
\begin{equation*}
u^{\sigma}=\sum_{k=1}^N u_k^\sigma, \quad \bbf^\sigma =\sum_{k=1}^N \bbf^\sigma_k,
\end{equation*}
where $\sigma$ represents the DOF in the FR context, i.e. the solution points to evaluate the polyomials.
In place of doing every calculation in each element $\Omega_k$ a reference element $\Omega_s$ is chosen.
Each element $\Omega_k$ is mapped to the reference element $\Omega_s$ and all calculation are done in $\Omega_s$.
The initial equation \eqref{eq:hyp_prob}
can be transformed to the following governing equation in the reference domain;
\begin{equation}\label{eq:hyp_prob_standard}
 \partial_t \tilde{u}+\div( \tilde{\bbf})=0 \Longleftrightarrow  \partial_t \tilde{u}+ \nabla_{rs}\cdot \tilde{\bbf}=0,
\end{equation}
where $\div$ is the divergence in the variables $r$, $s$  in the computational space. To clarify the notion we use  $\nabla_{rs}\cdot$
until the end of this section. 
$\tilde{u},\; \tilde{\bbf}$ are the transformed $u$ and $\bbf$.
The quantities $\tilde{u}$ and $\tilde{\bbf}$ can be directly calculated using the element mapping.
We suppress this dependence in the following to simplify the notation, see \cite{castonguay2012newclass} for details.\\
Let $P_p(\Omega_S)$ denote  the space of polynomials of degree less than $p$ on $\Omega_s$ and $R_p(\Gamma_S)$ be the polynomial space on the edges given by 
\begin{equation*}
 R_P(\Gamma_s)= \left\{\bbv \in L^2(\Gamma_s), \; \bbv|_{\Gamma_f}\in P_p(\Gamma_f), \;\forall \Gamma_f \right\}
\end{equation*}
where $\Gamma_f$ stand for the  edge $f$ of the reference element $\Omega_s$. 
The approximation $\tilde{u}^\sigma$ of the solution $u$ within the reference element $\Omega_s$ is done through a multi-dimensional polynomial of degree $p$, using
 the values of $\tilde{u}$ at $N_p=\frac{1}{2}(p+1)(p+2)$ solution points.
The solution approximation then reads;
\begin{equation*}
 \tilde{u}^\sigma(r,s,t)=\sum_{i=1}^{N_p} \tilde{u}^\sigma_i l_i(r,s)
\end{equation*}
where $\tilde{u}_i^\sigma$ is the value\footnote{We neglect the dependence on the mapping here again.} of $\tilde{u}$
at the solution point $i$ and $l_i(r,s)$ is the multidimensional Lagrange polynomial associated with the solution points $i$ in the
reference element $\Omega_s$.

 We now detail the main idea of FR. 
A simple  approach is to also approximate the flux function $f$ 
 by a polynomial $\tilde{\bbf}^\sigma=(\tilde{f}^\sigma, \tilde{g}^\sigma)$. To build this approximation, a first polynomial decomposition of $\tilde{f}^\sigma $ is set as
\begin{equation}\label{eq:flux_diverg}
 \tilde{f}^{\sigma,D}= \sum_{i=1}^{N_p} \tilde{f}^\sigma_i l_i(r,s), \quad  \tilde{g}^{\sigma,D}= \sum_{i=1}^{N_p} \tilde{g}^\sigma_i l_i(r,s) ,
\end{equation}
where the coefficients of these polynomials, respectively denoted by $\tilde{f}_i^\sigma$ and $\tilde{g}_i^\sigma$, are again evaluated at the solution points.
Since $\tilde{\bbf}^{\sigma, D}$ is always discontinuous at the boundary, it is called \emph{discontinuous flux}. 
To overcome / reduce that problem a further term $\tilde{\bbf}^{\sigma, C}$
is then added to $\tilde{\bbf}^{\sigma, D}$. $\tilde{\bbf}^{\sigma, C}$ is set to work directly on the boundaries of each element and corrects $\tilde{\bbf}^{\sigma, D}$ 
such that information of two neighboring
elements interacts and properties like conservation still hold in the discretisation.
We obtain for the approximation of the flux;
\begin{equation*}
 \tilde{\bbf}^\sigma= \tilde{\bbf}^{\sigma,D}+ \tilde{\bbf}^{\sigma,C}.
\end{equation*}
This gives \emph{Flux Reconstruction} its name. 
The selection/definition of these correction functions is essential. 
We now detail a possible construction for our special case. 
On each edge of the triangle,  a set of $N_{f_p}=p+1$ flux points are defined. These flux points are applied to couple 
the solution between neighboring elements. The correction function is then constructed as follows
\begin{equation}\label{eq:CORR_Fr}
 \tilde{\bbf}^{\sigma,C}=\sum_{f=1}^3 \sum_{j=1}^{N_{f_p}} \left[ \left( \tilde{\bbf} \cdot \bn \right)^{\sigma K}_{f,j}
 -\left(\tilde{\bbf}^{\sigma,D} \cdot \bn \right)_{f,j} \right] \bs h_{f,j}(r,s)
\end{equation}
The indices $f,j$ correspond to a quantity at the flux point $j$ of face $f$. Thus, in our case it is $1\leq f\leq 3$ and 
$1\leq j\leq N_{f_p}$. The term $\left(\tilde{\bbf}^{\sigma,D} \cdot \bn \right)_{f,j}$ is the normal component of the transformed discontinuous
flux at the flux point $j,j$, whereas $\left( \tilde{\bbf} \cdot \bn \right)^{\sigma K}_{f,j}$ is a normal transformed numerical flux 
computed at flux point $f,j$. We compute it by evaluating the multiply defined values of $u^\sigma$ at each flux point. More precisely,
we first define by $u^{\sigma,-}$ the value of $u^\sigma$ computed in the current element and by $u^{\sigma,+}$ its 
value computed using the information from the adjoint element that shares the same flux point $j$. This couples two neighboring elements and the information between them.
We then evaluate $u^{\sigma,+}$ and $u^{\sigma,-}$ at each flux point and compute 
$\left( \tilde{\bbf}(u^{\sigma,-},u^{\sigma,+})  \cdot \bn \right)^{\sigma K}_{f,j}$.
Finally, $\bs h_{f,j}(r,s)$ has to be explained. This is a vector correction function associated with the flux points $f,j$ 
and that lie in the Raviart-Thomas space $RT_p(\Omega_s)$ of order $p$. Thus $\bs h$ fulfills the following two properties 
\begin{equation}\label{eq:property_RT}
 \begin{aligned}
  &\nabla_{rs}\cdot \bs h_{f,j}  \in P_p(\Omega_s)\\
  &\bs h_{f,j}\cdot \bn \in R_p(\Gamma_S)\\
 \end{aligned}
\end{equation}
and has also to satisfy 
\begin{equation}\label{eq:Condition_COR_FR}
  \bs h_{f,j}(\bs r_{f_2,j_2}) \cdot \bn_{f_2,j_2}=\begin{cases}
                                                   1 \quad \text{ if } f=f_2 \text{ and } j=j_2,\\
                                                   0 \quad \text{ if } f\neq f_2 \text{ or } j\neq j_2.
                                                  \end{cases}
\end{equation}

Because of \eqref{eq:Condition_COR_FR} it follows 
that \begin{equation*}
      \tilde{\bbf}^{\sigma,C}(\bs r_{f,j})\cdot \bn_{f,j}= \left[  \left( \tilde{\bbf} \cdot \bn \right)^{\sigma K}_{f,j}
 -\left(\tilde{\bbf}^{\sigma,D} \cdot \bn \right)_{f,j}  \right]=:\alpha_{f,j}.
     \end{equation*}
     We also get 
$ \tilde{\bbf}^\sigma(\bs r_{f,j}) \cdot \bn_{f,j} = \left( \tilde{\bbf} \cdot \bn \right)^{\sigma K}_{f,j}$ at each flux point $f,j$.  
Combining our results, the approximate solution values to the problem \eqref{eq:hyp_prob}  can be updated at the solution points from 
\begin{align*}
\frac{\operatorname{d}u}{\operatorname{d}t }&=-\left(\nabla_{rs}\cdot \tilde{\bbf}^\sigma\right)\Big|_{\bs r_i}\\
&=-\left(\nabla_{rs}\cdot \tilde{\bbf}^{\sigma,D} \right)\Big|_{\bs r_i}-\left(\nabla_{rs}\cdot \tilde{\bbf}^{\sigma,C} \right)\Big|_{\bs r_i}\\
&=-\sum_{k=1}^{N_p} \tilde{f}^\sigma_k \frac{\partial l_k(r,s) }{\partial r} \Big|_{\bs r_i}
  -\sum_{k=1}^{N_p} \tilde{g}^\sigma_k \frac{\partial l_k(r,s)}{\partial s} \Big|_{\bs r_i}
  -\sum_{f=1}^3 \sum_{j=1}^{N_{f_p}} \left[ \left( \tilde{\bbf} \cdot \bn \right)^{\sigma K}_{f,j}
 -\left(\tilde{\bbf}^{\sigma,D} \cdot \bn \right)_{f,j} \right] \nabla_{rs}\bs h_{f,j}(\bs r)\\
 &=-\sum_{k=1}^{N_p} \tilde{f}^\sigma_k \frac{\partial l_k(r,s) }{\partial r} \Big|_{\bs r_i}
  -\sum_{k=1}^{N_p} \tilde{g}^\sigma_k \frac{\partial l_k(r,s)}{\partial s} \Big|_{\bs r_i}
  -\sum_{f=1}^3 \sum_{j=1}^{N_{f_p}} \alpha_{f,j} \nabla_{rs}\bs h_{f,j}(\bs r).
\end{align*}
Defining our FR scheme reduces to select the distributions 
of flux points and solutions points, as well as the form of our correction functions. 
The choice leads to several numerical methods with different properties. 
In \cite{castonguay2012newclass} special attention is paid on conservation and linear stability which 
restricts again the set of correction functions, but we do not go further into details here. \\
Finally, we want to mention that some of the most  famous schemes are embedded in this framework  by a right choice of correction 
functions and point distributions. 
To give a concrete example, the nodal Discontinuous Galerkin Spectral Element Method 
of Gassner et al. \cite{gassner2011comparison, gassner2013skew} can be named.

%% file: 2_FR_to_RD.tex
\section{Connection between Flux Reconstruction and Residual Distribution}\label{sec:FR_Basic}


Instead of using a variational or integral form like in DG, FR schemes are applied in
their discretisation of the differential form \eqref{eq:steady_problem} as it is described 
in subsection \ref{subsec:FR}.
The flux function is approximated by a polynomial of degree $k+1$ 
denoted by $\bbfh$. The discretisation of our underlying problem \eqref{eq:steady_problem} reads
\begin{equation}\label{eq:FR_Ba1}
\div (\bbfh +\bs \alpha \nabla \bs\psi)= 0,
\end{equation}
where  $\bs \alpha \nabla \bs\psi$  
is our correction function with the scaling term $\bs \alpha=\hbbf-\bbfh \cdot \bn$. We change here the notation from subsection \ref{subsec:FR}
on purpose to clarify that we are dealing now with the general case. We can translate it  into triangles
by setting $ \nabla \bs\psi:= \bs h$.
 We get

\begin{equation}\label{eq:FR_Ba2}
 \div \left(\bbfh + \left(\hbbf-\bbfh \cdot \bn \right) \nabla \bs\psi\right)= 0.
\end{equation}
We 
derive conditions on $\nabla \bs\psi $
so that this approach fits in the RD framework and that our methods have the
desirable properties of conservation and stability.
First, let us focus on FR schemes. The main idea of the
FR schemes is that the numerical flux at the boundaries will
be corrected by functions in such manner that information of 
two neighboring elements interacts and properties like conservation 
hold also in their discretisations. 
Let us  consider our discretisation \eqref{eq:FR_Ba2}. 
If we apply a Galerkin approach in every element $K$, then 
we obtain that for any  $\bs v^h \in V^h$ the relation
\begin{equation}\label{eq:FR_DG}
 \int_K \bs v^h \cdot \div \left(\bbfh  + \left(\hbbf (\bu^h, \bu^{h,-} )-\bbfh \cdot \bn \right) \nabla \bs\psi\right)
  \operatorname{d} \bxx =0
\end{equation}
has to be fulfilled. Using the Gauss theorem in the above equation yields
\begin{equation}\label{eq:Gauss_FR}
 -\int_K \bs \nabla v^h \cdot 
 \left(\bbfh + \bs \alpha \nabla \bs\psi\right) \operatorname{d} \bxx
 +\int_{\partial K} \bs v^h \cdot  \left( \bbfh \cdot \bn +
 \left(\hbbf-\bbfh \cdot \bn \right) \nabla \bs\psi \cdot \bn \right) \operatorname{d} \gamma =0.
\end{equation}
To guarantee conservation we demand that the flux over the element boundaries 
should be expressed only by the numerical flux of elements sharing this boundary.
Therefore, we require that
\begin{equation*}
 \left( \bbfh \cdot \bn +
 \left(\hbbf(\bu^h, \bu^{h,-} ) -
 \bbfh \cdot \bn \right) \nabla \bs\psi \cdot \bn\right)= \hbbf(\bu^h, \bu^{h,-} ),
\end{equation*}
which implies
\begin{equation}\label{eq:first_condition}
 \nabla \bs\psi \cdot \bn\equiv 1
\end{equation}
on the boundary. The relation \eqref{eq:first_condition} 
yields us a first property on our correction function
$\nabla \bs\psi$. 
\begin{remark}\label{Re:First_condition}
The condition \eqref{eq:first_condition} can be further weaken. Since we are using 
quadrature rules to evaluate the integrals \eqref{eq:FR_DG} or \eqref{eq:Gauss_FR}, \eqref{eq:first_condition}
has to be be only fulfilled at the quadrature points. 
We can guarantee the property \eqref{eq:first_condition} in
two dimension $(d=2)$ when using functions $\nabla\bs \psi$ lying 
in the lowest order Raviart-Thomas space \cite{raviart1977mixed}, 
up to some scaling.  The relation \eqref{eq:first_condition} is
then automatically fulfilled as a basic property of this function space.
In \cite{castonguay2012newclass} the authors already considered  
the Raviart-Thomas elements focusing on triangles.
We are considering the more general case of polygons. 
\end{remark}

To demonstrate the connection between FR and RD and to build an numerical scheme we have to 
apply quadrature formulas to evaluate the continious integrals. 

\subsection{From Flux Reconstruction to Residual Distribution Schemes}\label{sec:FR_to_RD}
In this part of the paper we show the connection between Flux Reconstruction and Residual Distribution Schemes. The key is a proper definition of the residuals.
If one has again a look on the formulation of the residuals  
\eqref{eq:DG_Variational}-\eqref{eq:DG_Residual_boundary} from section \ref{sec:RD_Basic} and compare them now with the formulations of 
\eqref{eq:FR_DG}- \eqref{eq:Gauss_FR}, it can be noticed that the equations  share  similar structures.
By passing from  integrals to quadrature formulas and splitting $\bs v^h$ along $\{\phi_\sigma\}_{\sigma\in S}$, we can define the residuals
in the following manner.
\begin{equation}\label{eq:DG_Residual_inner_FR}
\Phi_\sigma^{K,FR}(\bu^h) := \oint_K  \phi_\sigma 
 \cdot \div \left(\bbfh  + \left(\hbbf (\bu^h, \bu^{h,-} )-\bbfh \cdot \bn \right) \nabla \bs\psi\right) \operatorname{d} \bxx
\end{equation}
An other approach is to use Gauss formula (integration-by-parts/summation-by-parts),
leading to
\begin{equation}\label{eq:DG_Residual_inner_FR_2}
\begin{aligned}
 \Phi_\sigma^{K,FR}(\bu^h) :=& -\oint_K  \nabla  \phi_\sigma \cdot  \bbfh \operatorname{d} \bxx + 
\oint _{\partial K} \phi_\sigma \left(
 \bs \alpha  \nabla \bs\psi \cdot \bn 
+\bbfh \cdot \vecn \right) \operatorname{d}  \gamma
-\oint_K  \nabla \phi_\sigma \cdot \bs \alpha  \nabla \bs\psi \operatorname{d} \bxx \\
\stackrel{\eqref{eq:first_condition}}{=}&
-\oint_K \nabla \phi_\sigma \cdot \bbfh \operatorname{d} \bxx +  
\oint _{\partial K} \phi_\sigma \bs \hbbf (\bu^h, \bu^{h,-} ) \operatorname{d} \gamma
\underbrace{-\oint_K  \nabla \phi_\sigma \cdot \bs \alpha  \nabla \bs\psi
\operatorname{d} \bxx}_{:=r_\sigma}.
\end{aligned}
\end{equation}
Recalling property \eqref{eq:first_condition}, the boundary residuals reads
\begin{equation}\label{eq:DG_Residual_boundary_FR}
\Phi_\sigma^{K,FR}(\bu^h)=\oint_\Gamma \phi_\sigma  
 \left( \hbbf(\bu^h, \bu_b)-\bbf^h\cdot \vecn \right) \operatorname{d} \gamma.
\end{equation}
Comparing the residuals \eqref{eq:DG_Residual_inner_FR} and 
\eqref{eq:DG_Residual_inner_FR_2} with the residuals \eqref{eq:DG_Residual_inner} 
of the DG scheme\footnote{Here, we neglect the fact that in our description 
of  DG  we did not approximate the flux function by a polynomial.}, we can write 
\begin{equation}\label{eq:Residual_FR_DG}
 \Phi_\sigma^{K,FR}(\bu^h)= \Phi_\sigma^{K,DG}(\bu^h)+r_\sigma.
\end{equation}
Furthermore, the conservation relation \eqref{eq:residual_element} 
directly provides a second property on $\nabla\psi$, explicitly

\begin{equation}\label{eq:second_condition}
 \sum\limits_{\sigma\in K} r_\sigma=-\sum\limits_{\sigma\in K}\oint_K  \nabla \bs \phi_\sigma \cdot \bs \alpha  \nabla \bs\psi \operatorname{d} \bxx= 0.
\end{equation}
If we apply the residuals \eqref{eq:DG_Residual_inner_FR}-\eqref{eq:DG_Residual_boundary_FR} on our underlying steady state problem \eqref{eq:steady_problem}-\eqref{eq:initial condition}, 
we are able to write  our  model problem in the shape of \eqref{eq:Discretisatin_RD}.
For any $\sigma \in S$,  it reads
\begin{equation}\label{eq:Discretisatin_RD_FR}
 \sum\limits_{K\subset \Omega, \; \sigma \in K} \Phi^{K,FR}_\sigma (\bu^h) 
 +\sum\limits_{\ \Gamma \subset \partial \Omega, \; \sigma \in \Gamma} \Phi_\sigma^{\Gamma,FR} (\bu^h) =0. 
\end{equation}

With the definitions of the residuals \eqref{eq:DG_Residual_inner_FR}-\eqref{eq:Residual_FR_DG}
and the discretisation \eqref{eq:Discretisatin_RD_FR}, the Flux Reconstruction is
embedded within the RD framework. By ensuring that  conditions  \eqref{eq:first_condition} and \eqref{eq:second_condition} hold, the conservation relation \eqref{eq:residual_element} 
for our residuals is guaranteed and we are now able to use the theoretical results of  RD \cite{abgrall2006residual, abgrall2017some, abgrall2017some2, 
abgrall2017general, abgrall2017high} for the FR schemes under consideration. 
Naturally, the conservation properties of Flux Reconstruction schemes also hold and also stability results will 
transfer from the RD framework to FR.
\begin{remark}
If we are considering a two dimensional problem \eqref{eq:steady_problem},
our approach does not restrict the splitting of the domain $\Omega$ to a 
specific geometric structure like triangles or rectangles. The results
are valid more generally for all polygons.
This approach then extends the results of 
\cite{ castonguay2012newclass,huynh2014high,mengaldo2015dealiasing, vincent2011newclass} 
on FR to general grids. 
\end{remark}

Before we focus on entropy stability of our FR methods 
we shortly repeat some well-known results of RD schemes 
from  \cite{abgrall2017general, abgrall2017some, abgrall2017some2}
and references therein.

%% file: 3_Transformation.tex
\section{Transformation Results to Flux Reconstruction}\label{sec:Transformation}

As it is described inter alia in \cite{abgrall2017some}
for the RD schemes, a generalization of the classical Lax-Wendroff theorem is valued. It  transfers naturally to our FR formulation in RD \eqref{eq:DG_Residual_inner_FR}-\eqref{eq:Residual_FR_DG}.
\begin{theorem}[Theorem 2.2 of \cite{abgrall2003high}]\label{Theorem2_2}
 Assume the family of meshes $\TT =(\TT_h)$ is shape regular. We assume that the residuals  $\{\Phi_\sigma^{K}\}_{\sigma\in K}$ for $K$ an element or a boundary element of $\TT_h$ satisfy:
 \begin{itemize}
  \item For any $M\in \R^+$ , there exists a constant $C$ which depends only on the family of 
  meshes $\TT_h$ and $M$ such that for any $\bbu^h \in V^h$ with $||\bbu^h||_{\infty} \leq M$,
  then 
  \begin{equation*}
   ||\Phi_\sigma^{K}(\bbu^h|_{K})||\leq C \sum_{\sigma,\sigma'\in K}
   ||\bbu_\sigma^h-\bbu_{\sigma'}^h||.
  \end{equation*}
\item 
The conservation relations \eqref{eq:residual_element} and \eqref{eq:residual_boundary}
hold. 
 \end{itemize}
If there exists a Constant $C_{\max}$ such that the solutions of the scheme
\eqref{eq:Discretisatin_RD} (or \eqref{eq:Discretisatin_RD_FR}) satisfy $||\bbu^h||_\infty\leq 
C_{\max}$ and a function $\bv \in L^2(\Omega)^d $ such that $(\bbu^h)_h$ 
or at least a sub-sequence  converges to $\bv$ in $L^2(\Omega)^d$, then $\bv$ is a
weak solution of \eqref{eq:steady_problem}.
\end{theorem}
In view of the proof, the following relation is essential and 
can be derived from the conservation relations \eqref{eq:residual_element} 
and \eqref{eq:residual_boundary}. 
For any $\bv^h \in V^h$ written as 
$\bbv^h= \sum_{\sigma\in S} \bv_\sigma \phi_\sigma$:
\begin{equation}\label{eq:Conservation_relation_proof}
 \begin{aligned}
  0=& - \oint_{\Omega} \nabla \bv^h \cdot \bbf(\bbu^h) \operatorname{d} \bxx
  +\oint_{\partial \Omega}\bv^h\left( \hbbf(\bbu^h,\bbu_b)-\bbf (\bbu^h) \cdot \vecn \right) \operatorname{d} \gamma \\
  & +\sum_{\ee \in \EE^h} 
  \oint_{\ee} \{ \bv^h \} \hbbf(\bbu^h,\bbu^{h,-})\operatorname{d} \gamma
  +\sum_{K\in \Omega} \frac{1}{\# K} \left(\sum_{\sigma, \sigma' \in K} 
  \left(\bv_\sigma-\bv_{\sigma'} \right) \left(\Phi_\sigma^K(\bbu^h)-\Phi_{\sigma}^{K,DG} (\bbu^h)\right)\right)\\
  &+ \sum_{\Gamma \in  \partial 
  \Omega} \frac{1}{\# \Gamma} \left(\sum_{\sigma, \sigma' \in \Gamma} 
  \left(\bv_\sigma-\bv_{\sigma'} \right) \left(\Phi_\sigma^\Gamma(\bbu^h)-\Phi_{\sigma}^{\Gamma,DG} (\bbu^h)\right)\right),
 \end{aligned}
\end{equation}
with $\Phi_\sigma^{\bullet,DG}$ as in \eqref{eq:DG_Residual_inner} 
and \eqref{eq:DG_Residual_boundary}.

A consequence of \eqref{eq:Conservation_relation_proof}
is the following entropy inequality:
\begin{proposition}[Proposition 3.2 from \cite{abgrall2017some}]\label{Proposistion}
 Let $(U, \bbg)$ be an couple entropy-flux for \eqref{eq:steady_problem} and $\hbbg$ a numerical entropy flux consistent with $\bbg \cdot \vecn$. Assume that the residuals satisfy:
 \begin{equation}\label{conditon_Stability}
 \begin{aligned}
    \sum_{\sigma\in K} \est{\nabla_{\bbu} U(\bbu_\sigma),\Phi_\sigma^K} \geq& \oint_{\partial K} \hbbg(\bbu^h,\bbu^{h,-})
  \operatorname{d} \gamma &\text{ for any element K,}\\
   \sum_{\sigma\in \ee} \est{\nabla_{\bbu } U(\bbu_\sigma),\Phi_\sigma^\ee} \geq& \oint_{\partial \ee}
   \hbbg(\bbu^h,\bbu_b)-\bbg(\bbu^h)\cdot 
   \vecn )\operatorname{d} \gamma &\text{ for any boundary edge } \ee.
 \end{aligned}
 \end{equation}
Under the assumption of the theorem \ref{Theorem2_2}, the limit weak solution satisfies the following entropy inequality:
For any $\tau\in C^1(\Omega), \tau \geq 0,$
\begin{equation*}
 -\oint_\Omega \nabla \tau \cdot \bbg(\bbu)\operatorname{d} \bxx +\oint_{\partial \Omega^-}
 \tau \bbg(\bbu_b)\cdot \vecn \operatorname{d}  \gamma\leq 0.
\end{equation*}
\end{proposition}
The theorem  \ref{Theorem2_2} and the proposition \ref{Proposistion} 
ensure entropy stability. Therefore,
assumption \eqref{conditon_Stability} is essential. 
In \cite{abgrall2017general}, the property  \eqref{conditon_Stability} is further analysed and 
compared with the theory of Tadmor about entropy conservative/stable numerical flux functions 
\cite{tadmor1987numerical, tadmor2003entropy}. We can transfer this investigation, yielding a further condition on the correction functions. 
One can show that they have to satisfy an inequality in some sense to guarantee entropy stability.
However, this is more a theoretical condition and until now we do not see how this helps us to construct entropy conservative/stable FR schemes. 
For reasons of completeness and for future perspectives we develop this in the appendix. Another approach will be followed in the next section. 

\section{Entropy Conservative/Stable Flux Reconstruction Schemes} \label{sec:entropy}

We show in this section how to construct an entropy conservative scheme staring from our FR schemes 
$\Phi_\sigma^{K,FR}$ as the defined in \eqref{eq:DG_Residual_inner_FR_2} and \eqref{eq:Residual_FR_DG}.
From now on $\bbv$ represents the entropy variable 
$\nabla_\bbu U(\bbu)$.
Note that  since the entropy is strictly convex, the mapping $\bbu\to \bbv(\bbu)$ is one-to-one.

Here, we concentrate only on the the inner elements.
For a detailed analysis and the study of the boundary we strongly recommend \cite{abgrall2017some2}.
We  give further conditions on our correction functions in \eqref{eq:Residual_FR_DG} to get an entropy conservative/stable numerical FR scheme and presented a way to construct entropy stable FR schemes. This is done for the first time whereas \cite{castonguay2012newclass, vincent2011newclass, vincent2015extended} derive only linear stability.

Now, let $ \hbbg$ be the associated numerical entropy flux for the entropy variable $\bbv$ consistent 
with $\bs g \cdot \vecn$. 
As it is described in \cite{abgrall2017some2}
the entropy conditions can be formulated in the RD setting  for every inner element by the following formula.
\begin{equation}\label{eq:entropy_conservstion}
 \sum\limits_{\sigma \in K} \est{\bv_\sigma, \Phi^{CS}_\sigma} =\oint_{\partial K}\bs 
 \hbbg(\bv^h, \bv^{h,-}) \operatorname{d} \gamma
\end{equation}
To construct a new scheme $\Phi^{CS}_{\sigma}$
which fulfils  the condition \eqref{eq:entropy_conservstion},
we consider
\begin{equation}\label{eq:New_scheme}
 \Phi_\sigma^{CS} =
 \Phi_\sigma^{K,FR} + \tilde{\bs \tau}_\sigma\stackrel{\eqref{eq:Residual_FR_DG}}{=}
 \Phi_\sigma^{K,DG}+\bs r_\sigma + \tilde{\bs \tau}_\sigma,
\end{equation}
with some $ \tilde{\bs \tau}_\sigma$ that has to be built.
In order to guarantee the conservation property, we ask
\begin{equation}\label{eq:conservat_new}
 \sum\limits_{\sigma \in K} (\Phi^{K,FR}_\sigma+ \tilde{\bs \tau}_\sigma)
 =\oint_{\partial K} \bs \hbbf (\bv^h, \bv^{h,-}) \operatorname{d} \gamma,
\end{equation}
which implies 
\begin{equation}\label{eq:condition_r}
 \sum\limits_{\sigma\in K} \tilde{\bs \tau}_\sigma=0. 
\end{equation}
The question is now to construct $\tilde{\bs \tau}_\sigma$ under this constraint.
When  \eqref{eq:entropy_conservstion} holds, we have
\begin{equation}\label{eq:solution}
\sum\limits_{\sigma \in K} \est{\bv_\sigma, \tilde{\bs \tau}_\sigma}
=\oint_{\partial K} \bs \hbbg (\bv^h, \bv^{h,-}) \operatorname{d} \gamma- 
\sum\limits_{\sigma \in K} \est{\bv_\sigma, \Phi^{K,FR}_\sigma}:=\E.
\end{equation}
Thus, a solution to \eqref{eq:condition_r}-\eqref{eq:solution} is given by:
\begin{equation}\label{solution_r}
\tilde{\bs \tau}_\sigma=\alpha_\sigma (\bv_\sigma-\mean{\bv} ), \quad \alpha_\sigma
=\frac{\E}{\sum\limits_{\sigma \in K} (\bv_\sigma -\mean{ \bv} )^2 } ,  
\quad \mean{ \bv}= \frac{1}{\#K} \sum\limits_{\sigma \in K} \bv_{\sigma}.
\end{equation}
This can be seen by the following short calculation
\begin{align*}
  \sum\limits_{\sigma\in K} \tilde{\bs \tau}_\sigma&= \sum\limits_{\sigma\in K} \alpha_\sigma
  (\bv_\sigma- \mean{\bv})\stackrel{\eqref{solution_r}}{=}\sum\limits_{\sigma\in K} 
  \alpha_\sigma \bv_\sigma -\sum\limits_{\sigma \in K} \alpha_\sigma \bv_{\sigma}=0,\\
 E = \sum\limits_{\sigma \in K} \est{\bv_\sigma,\tilde{\bs \tau}_\sigma }
  &=\sum\limits_{\sigma \in K} \est{\bv_\sigma, \alpha_\sigma (\bv_\sigma-\mean{\bv} )}
  =  \sum\limits_{\sigma \in K} \alpha_\sigma 
  \left(\est{\bv_\sigma,  \bv_\sigma } -\est{\bv_\sigma, \mean{\bv}_\sigma}  \right)
  \\ &= \sum\limits_{\sigma \in K} \alpha_\sigma 
  \left(\est{\bv_\sigma,  \bv_\sigma } -2\est{\bv_\sigma, \mean{\bv}_\sigma}  +
  \est{ \mean{\bv}_\sigma, \mean{\bv}_\sigma}  \right)
 =\sum\limits_{\sigma \in K} \alpha_\sigma (\bv_\sigma -\mean{ \bv} )^2
\end{align*}
The scheme \eqref{eq:New_scheme} is entropy stable
by construction, but one can wonder about its accuracy.
We are using the second approach presented in  \cite{abgrall2017general},
since we can preserve the accuracy. 
%
%
The entropy flux $\bs g(\bu) \subset \R^d$ and the normal flux $\bbf$ fulfill
the following relation.
\footnote{Here, the expression $\est{\bv, \bbf}:= \bv^T \cdot \bbf$,
where the flux is interpreted as a $p\times d$ matrix
and the product is a $d$-row vector.}  
\begin{equation}\label{eq:entropy_potential}
 \bs g=\est{\bv, \bbf} -\bs \theta \text{ with } \nabla_{\bv} \bs \theta=\bbf
\end{equation}
The crucial point is the error $\E$ which has to be approximated as 
accurately as possible. For simplicity reason we are considering 
\begin{equation}\label{eq:error:FR}
 \tilde{\E}:= -\E= -\oint_{\partial K} \bs \hbbg (\bv^h, \bv^{h,-}) \operatorname{d} \gamma+ 
\sum\limits_{\sigma \in K} \est{\bv_\sigma, \Phi^{K,FR}_\sigma}=
-\oint_{\partial K} \bs \hbbg (\bv^h, \bv^{h,-}) \operatorname{d} \gamma+
\sum\limits_{\sigma \in K} \est{\bv_\sigma, \Phi^{K,DG}_\sigma+\bs r_\sigma}
\end{equation}
The entropy error $\tilde{\E}$ for DG was already investigated in \cite{abgrall2017general}.
We follow the steps which are analogous except for the additional terms $\bs r_\sigma$.
The numerical entropy $\hbbg$ will be defined later on. We hide the dependence of 
$\bbf^h $ on $\bv^h$ in the following to simplify the notation.
We further need the following two relations.

Using \eqref{eq:nabla_not} we have:
\begin{equation}\label{eq:chain_rule_system}
 \nabla \bv^h \cdot \bbfh = \div
 \left(\est{ \bv^h, \bbfh} \right) -\bv^h \cdot \nabla \bbfh,
\end{equation}
which can be also written in components.
It would read

\begin{equation*}
\sum\limits_{j=1}^d \sum\limits_{i=1}^p \left( \left(\frac{\partial }{\partial x_j}\bv^h_i \right)  f_{i,j} \right)
= \sum\limits_{j=1}^d \frac{\partial }{\partial x_j}  \left(\sum\limits_{i=1}^p (\bv^h_i f_{i,j}) \right) -
\sum\limits_{j=1}^d \sum\limits_{i=1}^p \left( \bv^h_i \left(\frac{\partial }{\partial x_j} f_{i,j}\right) \right).
\end{equation*}
We also get from \eqref{eq:Entrop} and \eqref{eq:Entrop_div}
\begin{equation}\label{eq_div_g}
 \bv(\bv^h) \cdot \nabla \bbf (\bv^h) 
 =\sum\limits_{j=1}^d \bv (\bv^h) \cdot \frac{\partial }{\partial x_j} \bbf_j(\bv^h)
 =\sum\limits_{j=1}^d \frac{\partial }{\partial x_j} g_j (\bv^h)=\div \bbg(\bv^h).
\end{equation}
Here, $\bv(\bv^h)$ is used to emphasis that condition \eqref{eq:Entrop} is only valued for the entropy variable 
$\bv$ and we can not assume directly that it holds for its interpolation. Therefore, 
we have to use the flux itself in this context and not the interpolated one.
Since we use the approximated flux function $\bbfh$ in our 
FR schemes we get a slight different investigation as in \cite{abgrall2017general}. 
However it does not change the major result. 
With \eqref{eq:chain_rule_system}, \eqref{eq_div_g} and  
Gauss's theorem\footnote{Here, we assume that the quadrature
has sufficient order or accuracy such that 
also the discrete version of the theorem is fulfilled. 
}
we are able to rewrite \eqref{eq:error:FR} as
\begin{align*}
 \tilde{\E} &\stackrel{\eqref{eq:chain_rule_system} }{ =}  
 - \oint_K \div \left(\est{ \bv^h, \bbfh} \right) \operatorname{d} \bxx + 
 \oint_K \bv^h \cdot \nabla \bbfh \operatorname{d} \bxx
 +\oint_{\partial K}  \left( \est{\bv^h, \hbbf(\bv^h, \bv^{h,-}) }  
 - \hbbg( \bv^h, \bv^{h,-} ) \right) \operatorname{d} \gamma +
 \sum_{\sigma\in K}\est{\bv_\sigma,\bs r_\sigma}\\
 &= 
 \underbrace{ -\oint_K \div \left(\est{ \bv^h, \bbfh} \right) \operatorname{d}
 \bxx + \int_K \div \left(\est{ \bv^h, \bbfh} \right) \operatorname{d} \bxx}_{:=SUR_1}
 -\int_K \div \left(\est{ \bv^h, \bbfh} \right) \operatorname{d} \bxx
\\
& \quad+ \underbrace{\oint_K  \left( \bv^h-\bv (\bv^h)\right)
\cdot \nabla \bbf(\bv^h) \operatorname{d} \bxx}_{:=SUR_2}
+ \underbrace{\oint_K  \bv^h
\cdot \nabla  \left(\bbfh- \bbf(\bv^h)\right) \operatorname{d} \bxx}_{:=SUR_3}
+ \oint_{  K} \bv (\bv^h) \cdot \nabla \bbf(\bv^h) \operatorname{d} \bxx 
\underbrace{-\oint_K\nabla \bv^h 
 \cdot \left(\bs\alpha \bs\psi \right)\operatorname{d} \bxx }_{:=CO}
\\
& \quad +\oint_{\partial K}  \left( \est{\bv^h, \hbbf(\bv^h, \bv^{h,-}) }  
 - \hbbg( \bv^h, \bv^{h,-} ) \right) \operatorname{d} \gamma 
 \\
&\stackrel{\text {Gauss} \& \eqref{eq_div_g} }{=} SUR_1+SUR_2+SUR_3+CO
\underbrace{-\int_{\partial K} \est{ \bv^h, \bbf(\bv^h)} \cdot 
\vecn \operatorname{d} \gamma +\oint_{\partial K} \est{ \bv^h, \bbf(\bv^h) }
\cdot \vecn \operatorname{d} \gamma }_{:=BO}\\
&\quad -\oint_{\partial K} \est{ \bv^h, \bbf(\bv^h)} \cdot \vecn \operatorname{d} \gamma 
+\oint_{\partial K} \bs g(\bv^h)\cdot \vecn \operatorname{d} \gamma  +
\oint_{\partial K}  \left( \est{\bv^h, \hbbf(\bv^h, \bv^{h,-}) }  
 - \hbbg( \bv^h, \bv^{h,-} ) \right) \operatorname{d} \gamma \\
&= SUR_1+SUR_2+CO+BO+  \oint_{\partial K} \left(  - \est{ \bv^h, \bbf(\bv^h)} \cdot \vecn  +  
\est{\bv^h, \hbbf(\bv^h, \bv^{h,-}) } + \bs g(\bv^h) \cdot \vecn - 
\hbbg (\bv^h, \bv^{h,-}) \right)\operatorname{d} \gamma
 \end{align*}
Assuming a smooth exact solution, we can use a quadrature 
formula of order $k$ for the three first volume terms
($SUR_1, \; SUR_2,\; SUR_3$) and obtain 
\begin{equation}\label{eq:Surf}
 SUR_1=O(h^{k+d+1}),\qquad SUR_2=  \Ol(h^{k+d+1}), \qquad SUR_3= \Ol(h^{k+d+1}).
\end{equation}
For the boundary term $BO$, we can get $ BO=\Ol(h^{k+d+1})$
using a quadrature formula of order $k+1$.
The last term has to be investigated. We use here for the numerical entropy flux 
\begin{equation}\label{eq:DG_g}
 \hbbg (\bv^h, \bv^{h,-}) = \est {\{ \bv^h \},\hbbf (\bv^h, \bv^{h,-}) }-
 \bs \theta (\{\bv^h \}).
\end{equation}
Applying \eqref{eq:entropy_potential} and \eqref{eq:DG_g} in the last term 
yields after some 
calculations to
\begin{equation*}
 - \est{ \bv^h, \bbf(\bv^h) } \cdot \vecn  +  
\est{\bv^h, \hbbf(\bv^h, \bv^{h,-}) } + \bs g(\bv^h) \cdot \vecn - \hbbg (\bv^h, \bv^{h,-} )
=\Ol(h^{2(k+1)})
\end{equation*}
and in total
\begin{equation*}
 \oint_{\partial K} \left(  - \est{ \bv^h, \bbf (\bv^h)} \cdot \vecn  +  
\est{\bv^h, \hbbf(\bv^h, \bv^{h,-}) } + \bs g(\bv^h) \cdot \vecn -
\hbbg (\bv^h, \bv^{h,-}) \right)\operatorname{d} \gamma=O(h^{d+2k+1})
\end{equation*}
with a suitable quadrature formula.

We now consider the last element term $CO$,
which has to approximate zero at least as $\Ol(h^{k+d+1})$. 
Indeed, if we assume that a quadrature of order $k$ for 
the volume integrals leads to $CO=\Ol(h^{d+k+1})$, then we can merge the errors and retrieve 
\begin{equation}\label{eq:total_error}
 \tilde{\E}= \sum\limits_{\sigma\in K} \est{\bv_\sigma, \Phi_\sigma} -\oint_{\partial K} 
\hbbg (\bv^h, \bv^{h,-} ) \operatorname{d} \gamma =\Ol(h^{d+k+1}),
\end{equation}
provided that a quadrature formula of order $k$ for the volume integrals 
and order $k+1$ for the boundary integrals are applied. 
This last equality furnishes a new condition for admissible 
correction functions and/or their quadrature rule.
This yields us to an extension of \textbf{Proposition 3.3.}
\begin{remark}\label{Re:Stability}
To obtain an entropy stable scheme, i.e. to have
\begin{equation}\label{eq:entropy_stable}
 \sum\limits_{\sigma \in K} \est{\bv_\sigma, \Phi^{ST}_\sigma} \geq 
 \oint_{\partial K}\bs  \hbbg(\bv^h, \bv^{h,-}) \operatorname{d} \gamma,
\end{equation}
we can combine the above approach and add an additional term to 
$\Phi^{CS}_\sigma$. More explicitly we can set: 
\begin{equation*}
 \Phi^{ST}_\sigma=\Phi^{CS}+\Psi_\sigma ,
\end{equation*}
where $ \Psi_\sigma$ satisfy $\sum\limits_{\sigma\in K}\Psi_\sigma=0 $ and $\sum\limits_{\sigma \in K} 
\est{\bv_\sigma, \Psi_\sigma}\geq0$. In \cite{abgrall2017general} two expressions for $\Psi$ can be found.
One contains jumps and the other one streamlines. 
\end{remark}
However, we are interested here in entropy conservative/stable FR schemes. Therefore, the inequality \eqref{eq:entropy_stable} 
should already hold for FR schemes of type
\eqref{eq:DG_Residual_inner_FR_2}.
Let  $\E^{DG}$ be the entropy error \eqref{eq:solution} of DG.
Up to assuming
\begin{equation}\label{theoretical_FR_COndition}
 \sum\limits_{\sigma\in K}\est{\bv_\sigma,\bs r_\sigma} \geq \E^{DG},
\end{equation}
we get automatically 
\begin{equation*}
 \sum\limits_{\sigma \in K} \est{\bv_\sigma, \Phi^{K,FR}_\sigma}
 -\oint_{\partial K} \bs \hbbg (\bv^h, \bv^{h,-}) \operatorname{d} \gamma
 =\sum\limits_{\sigma \in K} \est{\bv_\sigma, \bs r_\sigma}-
 \underbrace{\left(\oint_{\partial K} \bs \hbbg (\bv^h, \bv^{h,-}) \operatorname{d} \gamma
 -\sum\limits_{\sigma \in K} \est{\bv_\sigma, \Phi^{K,DG}_\sigma}\right)}_{\E^{DG}} \geq 0.
\end{equation*}

Thus, we are now able to present a way to build functions such 
that \eqref{eq:entropy_conservstion}  is naturally fulfilled.
To simplify our notation we include the scaling term $\bs \alpha$ in our correction functions 
$\nabla \bs \psi$ and use from now on $\nabla   \bs{\tilde{  \psi}} :=\bs \alpha \nabla \bs \psi$.
Therefore, conditions \eqref{eq:first_condition} and \eqref{eq:second_condition}
transfer to 
\begin{align}
 \nabla  \bs{\tilde{  \psi}}\cdot \bn &= \left(\hbbf(\bu^h, \bu^{h,-} ) -\bbfh \cdot \bn \right) \tag{\ref{eq:first_condition}*} \\
  \sum\limits_{\sigma\in K} r_\sigma &=-\sum\limits_{\sigma\in K}\oint_K  
 \nabla  \phi_\sigma \cdot \nabla  \bs{\tilde{  \psi}} \operatorname{d} \bxx= 0 \tag{\ref{eq:second_condition}*}.
\end{align}

For the construction of an entropy conservative scheme, 
we first add the term $\tilde{\bs \tau}$ to our FR schemes.
A solution of \eqref{solution_r} is then obtained. 
A way to determine a conservation preserving $\tilde{\bs \tau}$
is to read its formulation directly from the the ansatz \eqref{solution_r} 
used together with the above mentioned conditions, especially (\ref{eq:second_condition}*). 
On the one side we have
\begin{equation*}
 r_\sigma=\oint_K  
 \nabla  \phi_\sigma \cdot \nabla  \bs{\tilde{  \psi}} \operatorname{d} \bxx , 
\end{equation*}
but on the other side we also ask the entropy condition \eqref{eq:entropy_conservstion}. 
In order to fulfill this last condition, we may write $ r_\sigma $ as
\begin{equation*}
 r_\sigma =\alpha_\sigma( \E) \left(\bbv -\mean{\bv} \right),
\end{equation*}
where $\alpha_\sigma(\E)$ is defined as in \eqref{solution_r} and where we emphasized the dependence on the entropy error $\E$.
Bringing the two relations together we define our correction function by solving the following discrete Neumann problem.
\begin{equation}\label{eq:Neumann}
 \begin{aligned}
  \oint_K  
 \nabla  \phi_\sigma \cdot \nabla  \bs{\tilde{  \psi}} \operatorname{d} \bxx&= 
 \alpha_\sigma( \E) \left(\bbv -\mean{\bv} \right),\\
\nabla  \bs{\tilde{  \psi}}\cdot \bn &= \left(\hbbf(\bu^h, \bu^{h,-} ) -\bbfh \cdot \bn \right)
 \end{aligned}
\end{equation}
The equation \eqref{eq:entropy_conservstion} is then fulfilled 
and we naturally get entropy conservation for our FR schemes.
We can add jump or streamline terms to the schemes, 
as described in remark \ref{Re:Stability} or more explicitly 
in \cite{abgrall2017general}, to get entropy stability.

%% file: 4_Summary.tex
\section{Summary}

In this short paper we demonstrated the connection between Residual Distribution Schemes and the 
Flux Reconstruction approach. We saw
 that the FR schemes can be written as RD schemes. This link enables us to transform the well-known results about RD schemes into the FR framework. 
The crucial point is to derive suitable correction functions in FR.
From our previous analysis we can formulate the following main result. 
\begin{theorem}\label{Theorem_1}
We are considering the steady state problem \eqref{eq:steady_problem} together with the initial condition \eqref{eq:initial condition}. 
We are using a FR approach 
\begin{equation}\tag{\ref{eq:FR_Ba2}}
 \div \left(\bbfh + \bs \alpha \nabla \bs \psi\right)= 0.
\end{equation} 
with $\bs \alpha=\left(\hbbf-\bbfh \cdot \bn \right)$. Let us further assume that 
 our correction functions $\bs \alpha \nabla \bs \psi $ fulfill the following two conditions 
 \begin{itemize}
  \item \begin{equation}\tag{\ref{eq:first_condition}}
 \nabla \bs \psi \cdot \bn\equiv 1
\end{equation}
  \item 
  \begin{equation}\tag{\ref{eq:second_condition}}
 \sum\limits_{\sigma\in K} r_\sigma=-\sum\limits_{\sigma\in K}\oint_K  
 \nabla \bs \phi_\sigma \cdot \bs \alpha  \nabla \bs \psi \operatorname{d} \bxx= 0.
\end{equation}
 \end{itemize}
Then, our FR schemes are conservative and we are able to recast
them into the RD framework. All the results of section \ref{sec:Transformation} apply
automatically to our FR schemes.

If we further choose our correction functions so that 
\begin{equation}\tag{\ref{theoretical_FR_COndition}}
 \sum\limits_{\sigma\in K}\est{\bv_\sigma,\bs r_\sigma} \geq \E^{DG},
 \end{equation}
where $\E^{DG}$ is the entropy error \eqref{eq:solution} of DG,
then our FR scheme is additionally entropy stable.
\\
Finally, the correction functions are determined by solving the 
discrete Neuman problem \eqref{eq:Neumann}, leading to an entropy conservative FR scheme.
\end{theorem}
A naturally arising question is how to select this correction 
functions so that the conditions are fulfilled.
As explained in remark \ref{Re:First_condition},  
we search our correction functions $\bs \psi$ in the Raviart-Thomas space. 
Solving then the discrete Neumann problem \eqref{eq:Neumann} selects
entropy conservative FR schemes.

In \cite{vincent2015extended} 
Vincent et al. already described and developed FR schemes on triangular grids.
Their main problem was  to describe the correction functions. They  used
the Raviart-Thomas space and could  prove several conditions on their schemes
focusing on linear stability. 
They used in their analysis/construction a direct 
representation with flux points and solutions points on 
triangles whereas we apply an abstract approach.
Our advantage is that we do not use the geometrical structure of
the grid. Therefore, our results are valid for general polygons (in two space dimensions), 
and  include the schemes of \cite{vincent2015extended}.
However, this paper deals with the interpretation/transformation of FR schemes into the RD framework, allowing the use of theoretical results specific to RD in the context of FR. 
We are considering especially entropy stability and derive conditions
which are linked to the selection of correction functions. 
We further give an idea to build entropy conservative FR schemes. 
In a forthcoming paper we construct these FR schemes for general polygons and test them numerically against some benchmarks.

%% file: Appendix.tex
\section{Appendix}

\subsection{Entropy Stability- An Approach in the Sense of Tadmor}

In this appendix we  follow the steps from \cite{abgrall2017general} 
 and  have a closer look on the property \eqref{conditon_Stability} for our FR residuals 
 \eqref{eq:DG_Residual_inner_FR}-\eqref{eq:DG_Residual_boundary_FR}.
For simplicity reasons we assume in the following  that $K$ is a
fixed triangle. The results are nevertheless extendible to general
polytopes with degrees of freedom on the boundary of $K$.
We may consider a triangulation $\TT_K$ of $K$ elements whose
vertices are exactly the elements of $ S $.
We denote by $\hbbfe$ the flux between two DOFs $\sigma$ and $\sigma'$ and 
by $\vecn_{\sigma, \sigma_n}$ the normal vector on the direct edge between 
$\sigma$ and $\sigma_n$.

\begin{figure}[h]
\begin{center}
\subfigure[]{\includegraphics[width=0.45\textwidth]{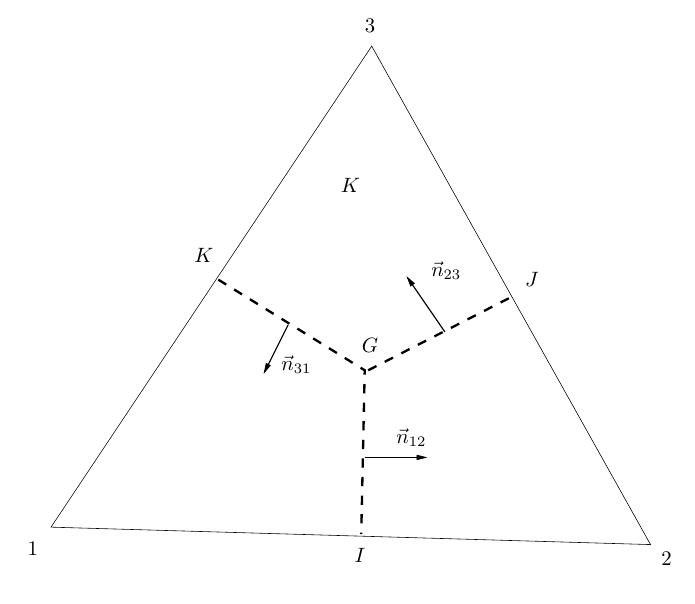}}
\subfigure[]{\includegraphics[width=0.45\textwidth]{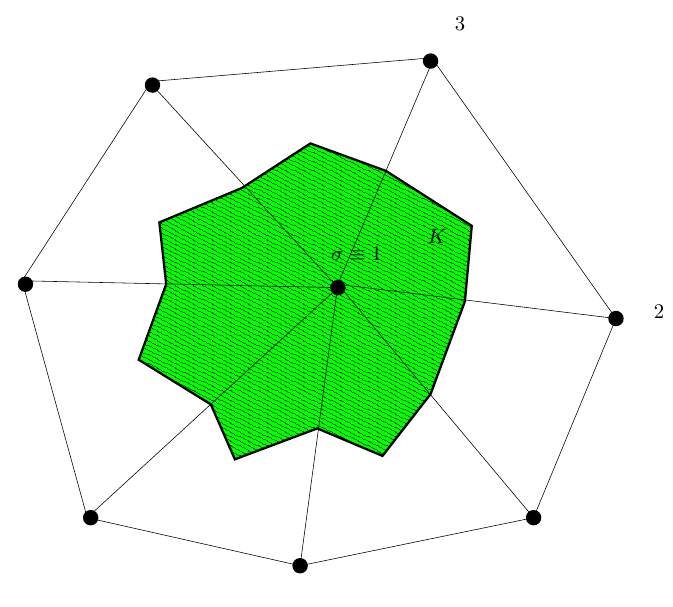}}
\end{center}
\caption{\label{fig:fv} Notations for the schemes. On the
left: Definition of the control volume for the degree of freedom $\sigma$.
  The vertex $\sigma$ plays the role of the vertex $1$ on the left
picture for the triangle K. The control volume $C_\sigma$ associated to
$\sigma=1$ is green on the right and corresponds to $1IGK$ on the left.
The vectors $\bn_{ij}$ are normal to the internal edges scaled by the
corresponding edge length.}
\end{figure}

In \cite{abgrall2017general} it is shown that we can split the residual in the following way;
\begin{equation}\label{eq:Residual_splitting}
 \Phi_\sigma^K= \sum_{\text{edges}[\sigma,\sigma']} \hbbfe+\hbbf_\sigma^b,
\end{equation}
with $\hbbf_\sigma^b= \oint_{\partial K} \phi_\sigma \hbbf_\sigma(u^h,u^{h,-})\operatorname{d} \gamma$.
Further properties and a detailed analysis can be found in \cite{abgrall2017general}.
As an additional example we derive the flux $\hbbfe$ for FR.
\begin{ex}[Flux Reconstruction schemes in the $\PP^1$ case]\label{ex_FR_Flux}  
The residuals are simply 
 \begin{equation*}
 \Phi_\sigma^{K,FR}(\bu^h) =-\oint_K \nabla \phi_\sigma \left( \bbfh + \bs \alpha  \nabla \bs \psi   \right)\operatorname{d} \bxx +  
\oint _{\partial K} \phi_\sigma \bs \hbbf (\bu^h, \bu^{h,-} ) \operatorname{d} \gamma
 \end{equation*}
The flux between two DOFs $\sigma$ and $\sigma'$ is given by
\begin{equation*}
 \hbbfe(\bbu^h,\bbu^{h,-})=   
\oint _{\partial K} ( \phi_\sigma -\phi_{\sigma'}) \hbbf (\bu^h, \bu^{h,-} ) \operatorname{d} \gamma
-\oint_K \nabla ( \phi_\sigma-\phi_{\sigma'}) \cdot \left( \bbfh + \bs \alpha  \nabla \bs \psi   \right)\operatorname{d} \bxx ,
\end{equation*}
where the equality
$ \nabla ( \phi_\sigma-\phi_{\sigma'})=\frac{\vecn_{\sigma \sigma'}}{|K|}$ 
comes from the  geometry (\red{See figures \ref{fig:fv} and \ref{control:DG}} ).
\begin{figure}[h]
\begin{center}
\includegraphics[width=0.5\textwidth]{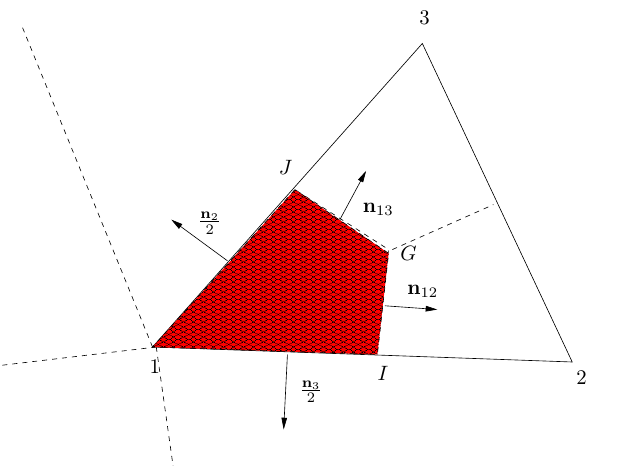}
\end{center}
\caption{\label{control:DG} Representation of the control volume
associated to DOF $1$.}
\end{figure}
We obtain 
\begin{equation}\label{eq:Flux_example}
 \hbbfe(\bbu^h,\bbu^{h,-})=   
\oint _{\partial K} ( \phi_\sigma -\phi_{\sigma'}) \hbbf (\bu^h, \bu^{h,-} ) \operatorname{d} \gamma
- \left(\oint_K  \left( \bbfh + \bs \alpha  \nabla \bs \psi   \right)\operatorname{d} \bxx \right)
\cdot \frac{\vecn_{\sigma \sigma'}}{|K|}.
\end{equation}
\end{ex}
\begin{remark}
 If the same quadrature formula is used
 on each edge, then $\oint_{\partial K} (\phi_\sigma-\phi_{\sigma'}) \operatorname{d} \gamma=0$ 
 and by introducing  the cell average $\overline{\bbu}^h$ of $\bbu^h$ and the
 average $\overline{\bbf}^h$ of the flux
 $\bbf^h$,  we can rewrite the flux as
 \begin{equation*}
 \hbbfe(\bbu^h,\bbu^{h,-})=   
\oint _{\partial K} ( \phi_\sigma -\phi_{\sigma'}) \left( 
\hbbf (\bu^h, \bu^{h,-} ) -\overline{\Bf}^h\cdot \vecn \right)
\operatorname{d} \gamma
- \left(\oint_K  \left( \bbfh + \bs \alpha  \nabla  \bs \psi   \right)\operatorname{d} \bxx \right)
\cdot \frac{\vecn_{\sigma \sigma'}}{|K|},
\end{equation*}
and the first term can be interpreted as a dissipation.
\end{remark} 
In Tadmors's work \cite{tadmor1987numerical, tadmor2003entropy} about 
entropy stability, sufficient conditions are derived for the schemes to
be entropy stable/conservative. 
In particular, it is shown that the numerical entropy flux needs to 
fulfill some dissipation inequality. We deduce here equivalent conditions for our setting, paying attention to the correction function that plays an important role. 
Using equation \eqref{eq:Residual_splitting} we can write our FR residuals as
\begin{equation}\label{eq:Residual_splitting_FR}
  \Phi_\sigma^{K,FR}= \sum_{\text{edges}[\sigma,\sigma']} \hbbfe^{FR}+\hbbf_\sigma^{b,FR},
\end{equation}
with $\hbbf_\sigma^{b,FR}=\oint_{\partial K} \phi_\sigma \hbbf(\bbu^h, \bbu^{h,-}) \operatorname{d} \bxx$.
From proposition \ref{Proposistion} we know that the condition
\eqref{conditon_Stability} is sufficient to guarantee entropy stability. Let us 
now insert \eqref{eq:Residual_splitting_FR} into it. 
We obtain\footnote{Here, we considered directly the entropy variable $\bbv$ similar to \cite{tadmor1987numerical} instead of $\bbu$.
Furthermore, we deal with an oriented graph. Given two vertices of this graph $\sigma$ and $\sigma'$, we write $\sigma> \sigma'$ a direct edge and we shorten the notation by $\sum\limits_{\sigma>\sigma'}:= \sum\limits_{\sigma \in K}\left(\sum\limits_{\sigma'\in K|\sigma>\sigma'} \right)$.}
\begin{equation}\label{eq:Asumption}
 \begin{aligned}
  \sum\limits_{\sigma\in K}\est{\bbv_\sigma, \Phi_\sigma^{K,FR}}=& 
 \sum\limits_{\sigma\in K}\est{\bbv_\sigma, \sum_{[\sigma,\sigma']} \hbbfe^{FR}+\hbbf_\sigma^{b,FR}}\\
 =& \sum\limits_{\sigma\in K} \est{\bbv_\sigma, \sum_{[\sigma,\sigma']} \hbbfe^{FR}}
  + \sum\limits_{\sigma\in K} \est{\bbv_{\sigma}, \oint_{\partial K} 
  \phi_\sigma \hbbf(\bbv^h, \bbv^{h,-}) \operatorname{d} \bxx} 
 \\
 =&  \frac12\sum_{\sigma>\sigma'} \est{\bbv_\sigma-\bbv_{\sigma'},\hbbfe^{FR} }    +\oint_{\partial K} \est{\bbv^h,
 \hbbf(\bbv^h,\bbv^{h,-})} \operatorname{d} \gamma 
 \end{aligned}
\end{equation}
We recall that we want to tune the correction function such that the obtained scheme 
is stable. It reduces here to ask 
\begin{equation}\label{eq:Assumption2}
 \sum\limits_{\sigma\in K}\est{\bbv_\sigma, \Phi_\sigma^{K,FR}}	{\geq} \oint_{\partial K}
	\hbbg(\bbv^h,\bbv^{h,-})\operatorname{d} \gamma.
\end{equation}
We now derive a condition that guarantees \eqref{eq:Assumption2} to hold.
Introducing the potential $\Theta^h_K$ in $K$ by 
\begin{equation*}
\Theta^h_K:=\sum_{\sigma\in K} \theta_\sigma \phi_\sigma \text{ with }
\theta_\sigma=\est{\bbv_\sigma,\Bf (\bbv_\sigma)}-\bbg(\bbv).
\end{equation*}
where we used the fact that the entropy flux $\bbg$ is related to the flux $\Bf$ by $\bbg=\est{\bbv,\Bf}-\Theta$, and defining the numerical entropy flux $\hbbg$ in
this spirit yields to the following expression.
\begin{equation*}
\hbbg(\bbv^h,\bbv^{h,-}):= \est{\{\bbv^h\}, \hbbf(\bbv^h,\bbv^{h,-} )} -\{\Theta^h_K\} \cdot \vecn
\end{equation*}
Together with \eqref{Eq:Jump_Average} we may rewrite \eqref{eq:Asumption} as:
\begin{align*}
  &\frac12\sum_{\sigma>\sigma'} \est{\bbv_\sigma-\bbv_{\sigma'},\hbbfe^{FR} }    
  +\oint_{\partial K} \est{\bbv^h,
 \hbbf(\bbv^h,\bbv^{h,-})} \operatorname{d} \gamma 
 \geq \oint_{\partial K}
\left(
\est{\{\bbv^h\}, \hbbf(\bbv^h,\bbv^{h,-} )} -\{\Theta^h_K\} \cdot \vecn \right)\operatorname{d} \gamma\\
\Longleftrightarrow& \frac12\sum_{\sigma>\sigma'}  \est{\bbv_\sigma-\bbv_{\sigma'},\hbbfe^{FR} } 
+\oint_{\partial K} \est{\bbv^h,
 \hbbf(\bbv^h,\bbv^{h,-})} \operatorname{d} \gamma  +\oint_{\partial K}\Theta^h_K \cdot \vecn \operatorname{d}
 \gamma 
 +\frac{1}{2} \oint_{\partial K} [\Theta^h_K]\cdot \vecn \operatorname{d} \gamma\\
 &\geq  \oint_{\partial K}
\est{\{\bbv^h\}, \hbbf(\bbv^h,\bbv^{h,-} )}\operatorname{d} \gamma \\
\Longleftrightarrow& 
 \frac12 \sum_{\sigma>\sigma'} \est{\bbv_\sigma-\bbv_{\sigma'},\hbbfe^{FR} } +
 \oint_{\partial K}\Theta^h_K \cdot \vecn \operatorname{d}
 \gamma 
 +\frac{1}{2} \oint_{\partial K} [\Theta^h_K]\cdot \vecn \operatorname{d} \gamma
 -\frac{1}{2} \oint_{\partial K} \est{[\bbv^h],
 \hbbf(\bbv^h,\bbv^{h,-})} \operatorname{d} \gamma 
 \geq 0
\end{align*}
Rearranging the last inequality we get:
\begin{equation}\label{eq:Con_entropy_stab}
 \underbrace{\frac12 \sum_{\sigma>\sigma'} \est{\bbv_\sigma-\bbv_{\sigma'},\hbbfe^{FR} } +
 \oint_{\partial K}\Theta^h_K \cdot \vecn \operatorname{d}
 \gamma}_{:=C_K}  -\underbrace{\frac{1}{2} \left( \oint_{\partial K} \left( \est{[\bbv^h],
 \hbbf(\bbv^h,\bbv^{h,-})} -[\Theta^h_K]  \right) \operatorname{d} \gamma   \right)}_{:=B_{\partial K}} \geq 0
\end{equation}
Taking now a closer look on our inequality \eqref{eq:Con_entropy_stab},
we may split this condition for entropy stability in two parts.
One part is actually working on the boundary of $K$, denoted
by $B_{\partial K} $, and one in $K$, denoted by $C_K$.

Let's pull a connection to the work of \cite{tadmor1987numerical}.
We recall that a numerical flux $\hbbf$ is entropy stable in the sense of Tadmor if
\begin{equation}\label{eq:Tadmor_entropy}
  \left( \est{[\bbv^h],
 \hbbf(\bbv^h,\bbv^{h,-})} -[\Theta^h_K]  \right) \leq 0
\end{equation}
holds. Thus combining $\eqref{eq:Tadmor_entropy} $ and 
the use of an entropy stable numerical flux in the sense
of Tadmor, we can guarantee that $B_{\partial K}\leq 0$.
We have only left to consider $C_K$. 
Let us study its first term. We have:

\begin{equation*}
 \frac12\sum_{\sigma>\sigma'} \est{\bbv_\sigma-\bbv_{\sigma'},\hbbfe^{FR} }
=\frac12 \sum_{\sigma>\sigma'}
\left( \est{\bbv_\sigma-\bbv_{\sigma'},\hbbfe^{FR} } -
\left( \theta _\sigma-\theta_{\sigma'} \right)\cdot \vecn_{\sigma,\sigma'}
\right)   +\frac{1}{2}\sum_{\sigma>\sigma'} \left( \theta _\sigma-\theta_{\sigma'} 
\right)\cdot \vecn_{\sigma,\sigma'}.
\end{equation*}
Furthermore, 
we have
\begin{equation}\label{eq:N}
 \frac{1}{2} \sum_{\sigma>\sigma'} 
 \left( \theta _\sigma-\theta_{\sigma'} \right)\cdot \vecn_{\sigma,\sigma'} =
 \sum_{\sigma\in K}  \theta_\sigma \cdot N_\sigma, 
\end{equation}
where the j-th component of $N$ writes $N_\sigma^j=\sum\limits_{edge [\sigma,\sigma']} \epsilon_{\sigma,\sigma'} 
\vecn_{\sigma,\sigma'}^j$ and where 
$\epsilon_{\sigma,\sigma'} $ is defined by 
\begin{equation*}
 \epsilon_{\sigma,\sigma'}:= \begin{cases}
                              0 \quad \text{ if $\sigma$ and $\sigma'$ are not on the same edge,}\\
                              1  \quad \text{ if $[\sigma,\sigma']$ is an edge and $\sigma\to \sigma'$ is direct,}      \\                     
                              -1 \quad \text{ if $[\sigma,\sigma']$ is an edge and $\sigma\to \sigma'$ is indirect.}
                             \end{cases}
\end{equation*}
Additionally, $N$ fulfills $\hbbf_\sigma^b=f(u^h)\cdot N_\sigma$ and we can write
 $N_\sigma:=-\oint_{\partial K} \phi_\sigma \cdot \vecn \operatorname{d} \gamma$. 
 Putting this in \eqref{eq:N} we get
 \begin{equation}\label{eq:N_2}
   \frac{1}{2}\sum_{\sigma\in K} \sum_{\sigma>\sigma'} 
 \left( \theta _\sigma-\theta_{\sigma'} \right)\cdot \vecn_{\sigma,\sigma'} =
 \sum_{\sigma\in K}  \theta_\sigma \cdot N_\sigma = -\sum_{\sigma\in K} \theta_\sigma \oint_{\partial K} \phi_\sigma \cdot \vecn \operatorname{d} \gamma
 = -\oint_{\partial K}\Theta^h_K \cdot \vecn \operatorname{d}
 \gamma
 .
 \end{equation}
By \eqref{eq:N_2} we get for $C_K$:
\begin{equation*}
  \frac12 \sum_{\sigma>\sigma'}
\left( \est{\bbv_\sigma-\bbv_{\sigma'},\hbbfe^{FR} } -
\left( \theta _\sigma-\theta_{\sigma'} \right)\cdot \vecn_{\sigma,\sigma'}
\right) +
 \oint_{\partial K}\Theta^h_K \cdot \vecn \operatorname{d}
 \gamma
  -\oint_{\partial K}\Theta^h_K \cdot \vecn \operatorname{d}
 \gamma \geq 0
\end{equation*}
Thus we finally get:
\begin{equation}\label{eq:Entropy_Stabiltiy_Tadmor_cond}
 \frac12 \sum_{\sigma>\sigma'}
\left( \est{\bbv_\sigma-\bbv_{\sigma'},\hbbfe^{FR} } -
\left( \theta _\sigma-\theta_{\sigma'} \right)\cdot \vecn_{\sigma,\sigma'}
\right)\geq 0.
\end{equation}
This condition is analogous to the one of Tadmor about entropy stability. 
Therefore, for an entropy stable flux function $\hbbfe^{FR}$, the entropy stability
is guaranteed.
However, note that $\hbbfe^{FR}$ depends on the correction functions.
We show it on our example \ref{ex_FR_Flux}.
\begin{ex}
 Putting \eqref{eq:Flux_example} into \eqref{eq:Entropy_Stabiltiy_Tadmor_cond}
 yields to the condition 
 \begin{align*}
&\frac12 \sum_{\sigma>\sigma'}
\left( \est{\bbv_\sigma-\bbv_{\sigma'},\hbbfe^{FR} } -
\left( \theta _\sigma-\theta_{\sigma'} \right)\cdot \vecn_{\sigma,\sigma'}
\right)\\
=&\frac12 \sum_{\sigma>\sigma'}
\left( \est{\bbv_\sigma-\bbv_{\sigma'},
\oint _{\partial K} ( \phi_\sigma -\phi_{\sigma'}) \hbbf (\bu^h, \bu^{h,-} ) \operatorname{d} \gamma
- \left(\oint_K  \left( \bbfh + \bs \alpha  \nabla \bs \psi   \right)\operatorname{d} \bxx \right)
\cdot \frac{\vecn_{\sigma \sigma'}}{|K|}} -
\left( \theta _\sigma-\theta_{\sigma'} \right)\cdot \vecn_{\sigma,\sigma'}
\right)\\
=& \frac12 \sum_{\sigma>\sigma'}
\left( \est{\bbv_\sigma-\bbv_{\sigma'},
\oint _{\partial K} ( \phi_\sigma -\phi_{\sigma'}) \hbbf (\bu^h, \bu^{h,-} ) \operatorname{d} \gamma
- \left(\oint_K  \left( \bbfh   \right)\operatorname{d} \bxx \right)
\cdot \frac{\vecn_{\sigma \sigma'}}{|K|}} -
\left( \theta _\sigma-\theta_{\sigma'} \right)\cdot \vecn_{\sigma,\sigma'}
\right)\\
&- \frac12 \sum_{\sigma>\sigma'} \est{\bbv_\sigma-\bbv_{\sigma'},
\left(\oint_K  \bs \alpha  \nabla \bs \psi \operatorname{d} \bxx \right)\cdot \frac{\vecn_{\sigma \sigma'}}{|K|}  }
{\geq}0
 \end{align*}
The first line of the last expression is again the same as for DG whereas the last line gives us special conditions on the correction functions which can be
tuned such that this inequality holds. 
Thus, the first step is the selection of  $\alpha \nabla \psi$ such that 
\begin{equation*}
 \frac12 \sum_{\sigma>\sigma'} \est{\bbv_\sigma-\bbv_{\sigma'},
\left(\oint_K  \bs \alpha  \nabla \bs \psi \operatorname{d} \bxx \right)\cdot \frac{\vecn_{\sigma \sigma'}}{|K|}  }\leq 0.
\end{equation*}
 \end{ex}
 In the above example we have seen that the correction function has a direct influence on entropy stability. We can now further restrict our correction functions so that the inequality \eqref{eq:Entropy_Stabiltiy_Tadmor_cond} holds.
A detailed analysis as well as more examples of those restrictions be found in
\cite{abgrall2017some2}.

\begin{remark}[Extension of Theorem \ref{Theorem_1}]
  Our investigation yields us an additional condition, so that we can extend our Theorem \ref{Theorem_1} by the the following:\\
     If we further choose our correction functions so that 
  $
\hbbfe^{FR}  $ is entropy stable in the sense of Tadmor \eqref{eq:Entropy_Stabiltiy_Tadmor_cond}
then our FR scheme is additionally entropy stable.
\end{remark}